\documentstyle[11pt]{article}

\input{amssym.def}
\input{amssym}

\newtheorem{theorem}{Theorem}

\newtheorem{definition}[theorem]{Definition}

\newtheorem{proposition}[theorem]{Proposition}
\newtheorem{remark}[theorem]{Remark}

\begin{document}

\def\cp{{\Bbb CP}^n}
\def\h{{\hbar}}
\def\gggg{{\bf g}}
\def\gh{\gggg_{\h}}
\def\slh{sl(2)_{\h}}
\def\glh{gl(2)_{\h}}
\def\gq{\gggg_{q}}
\def\gqq{\gggg_{q}^{\ot 2}}
\def\Ugh{U(gl(2)_{\h})}
\def\Ush{U(sl(2)_{\h})}
\def\Zgh{Z[U(gl(2)_{\h})]}
\def\Zsh{Z(U(sl(2)_{\h}))}
\def\uq{U_q(\gggg)}
\def\ugh{U(\gggg_{\h})}
\def\uqs{U_q(sl(n))}
\def\ug{U(\gggg)}
\def\C{{\Bbb C}}
\def\B{{\cal B}}
\def\cp{{\Bbb CP}}
\def\R{{\Bbb R}}
\def\Z{{\Bbb Z}}
\def\K{{\Bbb K}}
\def\ii{{\sqrt{-1}}}
\def\O{{\cal O}}
\def\A{{\cal A}}
\def\Ah{{\cal A}_\hbar}
\def\kO{{\K(\cal O)}}
\def\De{\Delta}
\def\de{\delta}
\def\h{{\hbar}}
\def\la{\lambda}
\def\La{\Lambda}
\def\co{{\cal O}}
\def\a{{\mu}}
\def\b{{\nu}}

\def\ll{L_{(2)}}
\def\oll{\overline L_{(2)}}
\def\MA{M_2{(\A)}}

\def\RE{RE(R)}
\def\sRE{sRE(R)}
\def\REh{RE(R)_{\h}}
\def\ip{I_+}
\def\im{I_-}
\def\ipm{I_{\pm}}
\def\ipgq{I_+(\gq)}
\def\imgq{I_-(\gq)}
\def\ipmgq{I_{\pm}(\gq)}
\def\ipL{I_+(\L)}
\def\imL{I_-(\L)}
\def\ipmL{I_{\pm}(\L)}
\def\Pp{P_+}
\def\Pm{P_-}
\def\Ppm{P_{\pm}}

\def\lp{\Lambda_+(V)}
\def\lm{\Lambda_-(V)}
\def\lpm{\Lambda_{\pm}(V)}
\def\lpl{\Lambda_+(\L)}
\def\lml{\Lambda_-(\L)}
\def\lpml{\Lambda_{\pm}(\L)}

\def\lpml{\Lambda_{\pm}^l(V)}
\def\lmp{\Lambda_{-}^p(V)}
\def\qm{q^{-1}}
\def\ve{\varepsilon}
\def\ep{\epsilon}
\def\al{{\alpha}}
\def\ah{{ A}_{\h}}
\def\aq{{ A}_{q}}
\def\ahq{{ A}_{\h,q}}
\def\ot{\otimes}
\def\om{\omega}
\def\Om{\Omega}
\def\Ml{M^{\lambda}}
\def\End{{\rm End\, }}
\def\Pic{{\rm Pic }}
\def\prePic{{\rm prePic }}
\def\Mod{{\rm Mod\, }}
\def\slpm{sl(p)-{\rm Mod}}
\def\slm{sl(n)-{\rm Mod}}
\def\Hom{{\rm End\, }}
\def\Sym{{\rm Sym\, }}
\def\Ker{{\rm Ker\, }}
\def\Im{{\rm Im\, }}
\def\Ob{{\rm Ob\, }}
\def\Gr{{\rm Gr\, }}
\def\SW{{SW(V)}}

\def\id{{\rm id\, }}
\def\dim{{\rm dim}}
\def\span{{\rm span\, }}
\def\tr{{\rm tr\, }}
\def\ad{{\rm ad\, }}
\def\cotr{{\rm cotr\, }}
\def\rk{{\rm rk\, }}
\def\Vect{{\rm Vect\,}}
\def\Fun{{\rm Fun\,}}
\def\vv{V^{\ot 2}}
\def\vvv{V^{\ot 3}}
\def\vl{V_{\la}}
\def\vm{V_{\mu}}
\def\lla{\bf{\la}}

\def\bea{\begin{eqnarray}}
\def\eea{\end{eqnarray}}
\def\be{\begin{equation}}
\def\ee{\end{equation}}
\def\nn{\nonumber}

\renewcommand{\theequation}{{\thesection}.{\arabic{equation}}}
\setcounter{equation}{0}

\title{Quantum line bundles on noncommutative sphere}
\author{D.~Gurevich\\
{\small\it ISTV, Universit\'e de Valenciennes, 59304 Valenciennes,
France}\\
\rule{0pt}{7mm}P.~Saponov,\\
{\small\it Theory Department of Institute for High Energy
Physics, 142284 Protvino, Russia}}

\maketitle

\begin{abstract}
Noncommutative (NC) sphere is introduced as a quotient  of the
enveloping algebra of the Lie algebra $su(2)$. Following
\cite{GS} and  using the Cayley-Hamilton identities we introduce
projective modules which are analogues of line bundles on the
usual sphere (we call them quantum line bundles) and define a
multiplicative structure in their family. Also, we compute a
pairing between quantum line bundles  and finite dimensional
representations of the NC sphere in the spirit of the NC index
theorem. A new approach to constructing the differential calculus
on a NC sphere is suggested. The approach makes use of the
projective modules in question and gives rise to a NC de Rham
complex being a deformation of the classical one.
\end{abstract}

{\bf AMS Mathematics Subject Classification, 1991 :} 17B37, 81R50


\section{Introduction}

One of basic notions of the usual (commutative) geometry is that
of the vector bundle on a variety. As was shown in \cite{Se}  the
category of vector bundles over a regular affine algebraic
variety $X$ is equivalent to the category of finitely generated
projective modules over the algebra ${\A}={\Bbb{K}(X)}$ which is
the coordinate ring of the given variety $X$ (a similar statement
for smooth compact varieties was shown in \cite{Sw}). Hereafter
$\Bbb{K}$ stands for the basic field (always $\Bbb{C}$ or
$\Bbb{R}$).

The language of projective modules is perfectly adapted to the
case of a noncommutative (NC) algebra $\cal{A}$. Any such (say,
right) $\cal{A}$-module can be identified with an idempotent
$e\in {M}_n(\cal{A})$ for some natural $n$. These idempotents
play the key role in all approaches to NC geometry, in
particular, in a NC version of the index formula of A.~Connes
(\cite{C}, \cite{L}).

The problem of constructing projective modules over physically
meaningful algebras is of great interest. The
$\Bbb{C}^*$-algebras (apart from commutative ones) are the mostly
studied from this viewpoint. Besides, there are known very few
examples. As an example let us evoke the paper \cite{Ri} where
projective modules over NC tori are studied (also cf. \cite{KS}
and the references therein). In the recent time a number of
papers have appeared dealing with some algebras (less standard
than NC tori and those arising from the Moyal product) for which
certain projective modules are constructed by hand (cf. for
example \cite{DL}, \cite{LM} and the references therein).

Nevertheless, there exists a natural  method suggested in
\cite{GS} of constructing projective modules over NC analogues of
$\Bbb{K}(\cal{O})$, where $\cal{O}$ is a generic\footnote{An
orbit is called generic if it contains a diagonal matrix with
pairwise distinct eigenvalues.} $SU(n)$-orbit in $su(n)^*$
including its "q-analogues" arising from the so-called reflection
equation (RE) algebra. This method is based on the
Cayley-Hamilton (CH) identity  for matrices with  entries
belonging to the NC algebras in question. (However, it seems very
plausible that other interesting examples of "NC varieties" can be
covered by this method. In particular, by making use of the CH
identity for super-matrices, cf. \cite{KT}, it is possible to
generalize our approach to certain super-varieties).

The idea of the method consists in the following. Consider a
matrix $L=\|l^i_j\|$ subject to the RE related to a Hecke
symmetry (cf. \cite{GPS})\footnote{A Hecke symmetry is a Hecke
type solution of the quantum Yang-Baxter equation. There exist
different types of Hecke symmetries: quasiclassical ones being
deformations of the classical flip and non-quasiclassical ones (a
big family of such Hecke symmetries was introduced in \cite{G}).
Let us note that a version of the CH identity  exists for any of
them independently on the type.}. Then it satisfies a polynomial
identity \be L^p+\sum_{i=1}^p a_i(L) L^{p-i}=0 \label{CHid} \ee
where coefficients $a_i(L)$ belong to the center of the RE
algebra generated by $l^i_j$.

Passing to a specific limit in RE algebra we get a version of the
CH identity for the matrix whose entries $l^i_j$ commute as
follows: \be [l^i_j,\,l^k_l]=\hbar(\delta_j^k\,l_l^i -
\delta^i_l\,l^k_j), \label{gln} \ee i.e. $l^i_j$ generate the Lie
algebra $gl(n)_\hbar$ where $\gh$ stands for a Lie algebra whose
Lie bracket equals $\h [\,\,,\,\,]$, where $[\,\,,\,\,]$ is the
bracket of a given Lie algebra $\gggg$. Introducing the parameter
$\hbar$ allows us to consider the enveloping algebras as
deformations of commutative ones. (Let us note that this type of
the CH identity was known since 80's, cf. \cite{Go}).

Similar to the general case the coefficients of the corresponding
CH identity belong to the center $Z[U(gl(n)_\hbar)]$ of the
enveloping algebra $U(gl(n)_\hbar)$. Therefore, by passing to a
quotient\footnote{Hereafter $\{X\}$ stands for the ideal
generated by a set $X$ in the algebra in question.}
$$
U(gl(n))/\{z-\chi(z)\}
$$
where $z\in Z[U(sl(n))]$ and
$$
\chi: Z[U(gl(n))]\to \Bbb{K}
$$
is a character we get a CH identity with numerical coefficients
\be L^p+\sum_{i=1}^p \alpha_i L^{p-i}=0,\quad
\alpha_i=\chi(a_i(L)).\label{numCH} \ee Assuming the roots of
this equation to be distinct pairwise we can assign an idempotent
(or what is the same one-sided projective module) to each root.
For quasiclassical Hecke symmetries (see footnote 2) these
projective modules are deformations of line bundles on the
corresponding classical variety. For this reason we call them
{\em quantum line bundles} (q.l.b.).

In this paper we constrain ourselves to the case arising from the
Lie algebra $gl(2)$. Namely, we describe a family of projective
modules over the algebra \be {\cal
A}_\hbar=U(sl(2)_{\hbar})/\{\Delta-\alpha\}, \label{quot} \ee
where  $\Delta$ stands for the Casimir element in the algebra
$U(sl(2)_{\hbar})$\footnote{Note that in fact we deal with the
Lie algebra $sl(2)_\hbar$ instead of $gl(2)_\hbar$ since the
trace of the matrix $L$ is always assumed to vanish. The studies
of such quotients were initiated by J.Dixmier (cf. \cite{H} and
the references therein). In certain papers compact form of this
algebra is called {\em fuzzy} sphere.}. We consider this algebra
as a NC counterpart of a hyperboloid.

In order to get a NC counterpart of the sphere we should pass to
the compact form of the algebra in question. However, it does not
affect the CH identity since it is indifferent to a concrete form
(compact or not) of the algebra.

To describe our method in more detail we begin with the classical
(commutative) case. Put the matrix
$$L=\left(\matrix{ix&-iy+z\cr -iy-z&-ix}\right)$$
in correspondence  to a point $(x,y,z)\in S^2$. This matrix
satisfies the CH identity (\ref{numCH}) where $p=2$,
$\alpha_1=0$, $\alpha_2=x^2+y^2+z^2={\rm const}\not=0$.

Let $\lambda_1$ and $\lambda_2=-\lambda_1$ be the  roots of
equation (\ref{numCH}). To each point $(x,y,z)\in S^2$ we assign
the eigenspace of the above matrix $L$ corresponding to the
eigenvalue $\lambda_l$, $l=1,2$. Thus, we come to a line bundle
$E_l$ which will be called {\em basic}.

Various tensor products of the basic bundles $E_l$, $l=1,2$, give
rise to a family of {\em derived or higher} line bundles
$$
E^{k_1,k_2}=E_1^{\ot k_1}\ot E_2^{\ot k_2},\,\,k_1,k_2=0,1,...
$$
Note, that certain line bundles of this family are isomorphic to
each other. In particular, we have  $E_1\ot E_2=E^{0,0}$ where
$E^{0,0}$ stands for the trivial line bundle. In general, the
line bundles
$$
E^{k_1,k_2} \quad {\rm and}\quad E^{k_1+l,k_2+l},\,\,l=1,2,...
$$
are isomorphic to each other. So, any line bundle is isomorphic
either to $E_1^{k}$ or to $E_2^{k}$ for some $k=0,1,...$ (we
assume that $E_1^{0}=E_2^{0}=E^{0,0}$). Finally, we conclude that
the Picard group $\Pic (S^2)$ of the sphere (which is the set of
classes of isomorphic lines bundles equipped with the tensor
product) is nothing but $\Bbb{Z}{}$ since any line bundle can be
represented as $E_1^{ k}$ with a proper $k\in \Bbb{Z}$, where we
put $E_1^{ k}=E_1^{\ot k}$ for $k>0$  and  $E_1^{ k}=E_2^{\ot
(-k)}$ for $k<0$.

It is worth emphasizing that we deal with an algebraic setting:
the sphere and total spaces of all bundles in question are
treated as real or complex affine algebraic varieties (depending
on the basic field $\Bbb K$).

Now, let us pass to the NC case. Using the above NC version of
the CH identity one can define NC analogues $E_l(\h)$, $l=1,2$,
of the line bundles $E_l$. We will call them {\em basic} q.l.b.
(they are defined in section 2).

Unfortunately, for a NC algebra $\A$ any tensor product of two or
more one-sided (say, right) $\A$-modules is not well-defined. So,
the construction of derived line bundles cannot be generalized to
a NC case in a straightforward way. Nevertheless, using the CH
identities for some "extensions" $L_{(k)},\,k=2,3...$ of the
matrix $L$ to higher spins (in the sequel $k=2\times$spin) we
directly construct NC counterparts $E^{k_1,k_2}(\h)$ of the above
derived line bundles. Thus, for any $k=2,3,...$ (with
$k_1+k_2=k$) we have $k+1$ {\em derived} q.l.b.

In section 3 we explicitly calculate the CH identity for the
matrix $L_{(2)}$ and state that such an identity exists for any
matrix $L_{(k)},\,\,k>2$.

However, as usual we are interested in projective modules (in
particular, q.l.b.) modulo natural isomorphisms. Thus, we show
that the q.l.b. $E^{1,1}(\h)$ is isomorphic to the trivial one
$E^{0,0}(\h)$ which is nothing but the algebra $\Ah$ itself. We
also conjecture that the q.l.b. \be E^{k_1,k_2}(\h)\quad {\rm
and}\quad E^{k_1+l,k_2+l}(\h)\label{isom} \ee are isomorphic to
each other. If it is so, any q.l.b. is isomorphic to
$E^{k,0}(\h)$ or $E^{0,k}(\h)$ similar to the commutative case.

Then we define an associative product in the family of q.l.b.
over the NC sphere in a natural way. By definition, the product
of two (or more) q.l.b. over NC sphere are the NC analogue of the
product of their classical counterparts. Otherwise stated, we set
by definition
$$
E^{k_1,k_2}(\h)\cdot E^{l_1,l_2}(\h)=E^{k_1+l_1,k_2+l_2}(\h)
$$
(in particular, we have $E_1(\h)\cdot E_2(\h)= E^{1,1}(\h)$ and
in virtue of the above mentioned result this product is
isomorphic to $E^{0,0}(\h)$). The family of all modules
$E^{k_1,k_2}(\h)$ equipped with this product is a semigroup. It
is denoted $\prePic(\Ah)$ and called prePicard. Assuming the
above conjecture to be true we get the Picard group $\Pic(\Ah)$
of the NC sphere (it is a group since any its element becomes
invertible)\footnote{Let us note that $K_0$ of the NC sphere
equipped only with the additive structure was calculated in
\cite{H}. We are rather interested in quantum line bundles. Once
they are defined other projective modules can be introduced as
their direct sums.}. All these notions are introduced in the
section 4. Moreover, in this section we compute the pairing
between  q.l.b. in question and irreducible representations of
the algebra $\Ah$ in the spirit of the NC index theorem.

In the last section we suggest a method  of constructing a
version of differential calculus on the NC sphere similar to that
on the commutative sphere. The method makes use of the above
projective modules instead of the Leibniz rule which is
habitually employed in this area. In contrast with the usual
approach giving rise to much bigger differential algebra than the
classical one our approach leads to the de Rham complex whose
terms are projective modules being deformations of the classical
ones.

To complete the Introduction we emphasize that our approach is in
principle applicable to NC analogue of any generic orbit in
$su(n)^*$ (and even to their "q-analogues" arising from the RE
algebra) but the calculation of the higher CH identities becomes
much more difficult as the degree of the CH identity for the
initial matrix becomes greater than 2. Actually, a paper
\cite{GLS} devoted to mentioned "q-analogues"  is in progress.

\medskip

{\bf Acknowledgment} The article was started during the authors'
visit to Max-Planck-Institut f\"ur Mathematik in Bonn. The
authors would like to thank the staff of the MPIM for warm
hospitality and stimulating atmosphere.

\section{Basic quantum line bundles on the NC sphere}
\setcounter{equation}{0}

Consider the Lie algebra $\glh$ which is generated by the elements
$a,b,c$ and $d$ satisfying the following commutation relations
$$
[a,b]=\hbar b,\quad [a,c]=-\hbar c,\quad  [a,d]=0,\quad
[b,c]=\hbar (a-d), \quad [b,d]=\hbar b,\quad [c,d]=-\hbar c.
$$
The nonzero numerical parameter $\hbar$ is introduced into the Lie
brackets for the future convenience. This parameter can be
evidently equated to one by the renormalization of the
generators. In this case we come to the conventional Lie algebra
$gl(2)$.

Let us form a $2\times 2$ matrix L whose entries are the above
generators
$$
L=\left(\matrix{a&b\cr c&d}\right).
$$
It is a matter of straightforward checking that this matrix
satisfies the following second order polynomial identity \be
L^2-(\tr +\hbar)L+(\triangle+\hbar\,\tr/2)\id=0 \label{CH} \ee
where
$$
\tr=\tr L=a+d,\quad   \triangle=ad-(bc+cb)/2.
$$
As $\hbar\rightarrow 0$ we get the classical CH identity for a
matrix with commutative entries.

In order to avoid any confusion we want to stress that we  only  deal with
 the enveloping algebra $U(\glh)$ (and some its quotients) and we
disregard its (restricted) dual object --- the algebra of
functions on the Lie group $GL(2)$. Our immediate aim
consists in constructing some "derived" matrices with
entries from $U(\glh)$ satisfying some "higher" CH
identities. It will be done by some sort of  coproduct
applied to the matrix $L$  and restricted onto the
symmetric component. However, we do  not use any
coalgebraic (and hence any Hopf) structure of the algebra
$U(\glh)$ itself.

\begin{remark}
A version of the CH identity for matrices with entries from
$U(gl(n))$ is known for a long time (cf. \cite{Go}). However,
traditionally one deals with the CH identity in a concrete
representation of the algebra $gl(n)$ while we prefer to work
with the above universal form of the CH identity. A way of
obtaining the CH identity by means of the so-called Yangians was
suggested in \cite{NT}. In \cite{G-T} another NC version of the
CH identity was presented. The coefficients of the polynomial
relation suggested there are diagonal (not scalar) matrices. But
such a form of the CH identity is not suitable for our aims.
\end{remark}

Since the elements tr and $\triangle$ from (\ref{CH}) belong to
the center $\Zgh$ of the algebra $\Ugh$ one can consider the
quotient
$$
\Ah=\Ugh/\{\tr, \,\triangle-\alpha\},\,\, \alpha\in \Bbb{K}.
$$
Taking in consideration the fact that the trace of the matrix $L$
vanishes we can also treat this algebra as quotient (\ref{quot}).
In what follows the algebra $\Ah$ will be called a NC variety (or
more precisely, a NC hyperboloid).

Being restricted to the algebra $\Ah$ the CH identity becomes a
polynomial relation in $L$ with numerical coefficients \be
L^2-\hbar L+\alpha\id=0. \label{numer} \ee Denote $\lambda_1$ and
$\lambda_2$ the roots of this equation that is
$$
\lambda_1=(\h-\sqrt{\h^2-4\alpha})/2,\quad \lambda_2=
(\h+\sqrt{\h^2-4\alpha})/2.
$$

Let us suppose that $\lambda_1\neq\lambda_2$. If $\hbar=0$ this
condition means that the cone corresponding to the case
$\triangle=\alpha=0$ is forbidden. However, if $\hbar\not=0$ we
have $\hbar^2-4\alpha\not=0$.

Abusing the language we say that $\lambda_1$ and $\lambda_2$ are
eigenvalues of the matrix $L$. Of course, this does not mean the
existence of an invertible matrix $A\in M_2(\Ah)$ such that the
matrix $A\cdot L\cdot A^{-1}$ becomes diagonal:
diag$(\lambda_1,\lambda_2)$.

\begin{remark} If an algebra $\A$ is not
a field, the invertibility of a matrix $A\in M_n(\A)$ is an
exceptional situation. As follows from the CH identity for the
matrix $L$ it is  invertible. However, it is not so even for
small deformations of $L$.
\end{remark}

It is easy to see that the matrices \be
e_{10}=(\lambda_2\id-L)/(\lambda_2-\lambda_1),\;\;
e_{01}=(\lambda_1\id-L)/(\lambda_1-\lambda_2)\in
M_2(\Ah)\label{ide} \ee are idempotents and $e_{10}\cdot
e_{01}=0$.

Denote $E_1(\h)$ and $E_2(\h)$ the projective modules (also
called q.l.b.) corresponding respectively to the idempotents
$e_{10}$ and $e_{01}$ in (\ref{ide}). Let us explicitly describe
these modules.

Let $V_{(k)}$ be the $k$-th homogeneous component $\Sym^k(V)$ of
the symmetric algebra $\Sym(V)$ of the space $V$. Thus,
$k=2\times$spin and $\dim(V_{(k)})=k+1$. Consider the tensor
product $V\ot \Ah$. It is nothing but the free right
$\Ah$-module  ${\Ah}^{\oplus 2}$. We can imagine the matrix $L$
(as well as any polynomial in it) as an operator acting from $V$
to $V\ot \Ah$ which can be presented in a basis $(v_1, v_2)$ as
follows:
$$
(v_1, v_2)\mapsto (v_1\ot a+ v_2\ot c, \,v_1\ot b+ v_2\ot d).
$$

Following  \cite{GS} we define the projective module
$E_l(\h),\,\,l=1,\,2$ as quotient of $V\ot \Ah$ over its
submodule generated by the elements \be v_1\ot a+ v_2\ot c-v_1\ot
\lambda_l, \quad v_1\ot b+ v_2\ot d-v_2\ot \lambda_l. \label{equa}
\ee Also,  the module $E_1(\h)$ (resp. $E_2(\h)$) can be
identified with the image of the idempotent $e_{10}$ (resp
$e_{01}$) which consists of the elements
$$
(v_1\ot a+ v_2\ot c-v_1\ot \lambda_m) f_1+(v_1\ot b+ v_2\ot
d-v_2\ot \lambda_m) f_2,\quad \forall\,f_1,f_2\in\Ah
$$
where $m=2$ for  $E_1(\h)$  and $m=1$ for $E_2(\h)$. (In a
similar way we can associate  left projective modules to the
idempotent in question.)

Now, consider the compact form of the NC variety in question.
Changing the basis
$$
x=i(d-a)/2=-ia,\,\,y=i(b+c)/2,\,\,z=(b-c)/2
$$
we get the following commutation relations between the new
generators \be [x,y]=\hbar z,\,\,[y,z]=\hbar x,\,\,[z,x]=\hbar y
\label{comp} \ee and the defining equation of the NC variety
reads now
$$
\triangle=(x^2+y^2+z^2)=\alpha.
$$

Thus, assuming $\K=\R$ and $\alpha>0$ we get a NC analogue of the
sphere, namely, the algebra
$$
\Ah=U(su(2)_\hbar)/\{x^2+y^2+z^2-\alpha\}
$$
However, the eigenvalues of equation (\ref{numer}) are imaginary
(for positive $\alpha$ and real and small enough $\h$) and as
usual we should consider the idempotents and related projective
modules over the field $\C$.

Completing this section we want to stress that equations
(\ref{equa}) are covariant w.r.t. the action of the group $G$
where $G=SU(2)$ or $G=SL(2)$ depending on the form (compact or
not) we are dealing with.

\section{Derived quantum line bundles}
\setcounter{equation}{0}

In this section we discuss the problem of extension of the matrix
$L=L_{(1)}$ to the higher spins and suggest a method of finding
the corresponding CH identities.

First, consider the commutative case. Let \be
\Delta(L)=L\ot\id+\id\ot L\in M_4(\Ah) \label{first} \ee be the
first extension of the matrix $L$ to the space $\vv$. If
$\lambda_1$ and $\lambda_2$ are (distinct) eigenvalues of $L$ and
$u_1, u_2\in V$ are corresponding eigenvectors, then the spectrum
of $\Delta(L)$ is
$$
2\lambda_1\; (u_1\ot u_1),\quad 2\lambda_2\; (u_2\ot u_2),\quad
\lambda_1 +\lambda_2\; (u_1\ot u_2 \;\;{\rm and }\;\; u_2\ot u_1)
$$
where in brackets we indicate the corresponding eigenvectors.

The commutativity of entries of the matrix $L$ can be expressed
by the relation
$$
L_1\cdot L_2=L_2\cdot L_1\quad{\rm where}\quad
  L_1=L\ot\id,\, L_2=\id\ot L.
$$

Rewriting (\ref{first}) in the form $\Delta(L)=L_1+L_2$ and
taking into account the CH identity for $L$
$$
0 = (L-\lambda_1\,\id)(L-\lambda_2\,\id) = L^2 -\a L
+\b\,\id,\quad \a=\lambda_1+\lambda_2,\;\b=\lambda_1\lambda_2
$$
we find
\begin{eqnarray*}
&&\Delta(L)^2 = \a\Delta(L)+2L_1\cdot L_2 - 2\b\,\id\\
&&\Delta(L)^3 = (\a^2 - 4\b)\Delta(L)+6\a L_1\cdot L_2 -
2\a\b\,\id.
\end{eqnarray*}
Note that $\a=0$ if $L\in sl(2)$. Upon excluding $L_1\cdot L_2$
from the above equations we get
$$
\Delta(L)^3-3\a\Delta(L)^2 +(2\a^2+4\b)\Delta(L) - 4\a\b\, \id=0.
$$
Substituting the values of $\a$ and $\b$ we can present this
relation as follows
$$
(\Delta(L) - 2\lambda_1)(\Delta(L) - \lambda_1 -\lambda_2)
(\Delta(L) - 2\lambda_2) = 0.
$$
Thus, the minimal polynomial for $\Delta(L)$ is of the degree 3.

A similar statement is valid for further extensions of the matrix
$L$:
$$
\Delta^2(L)=L\ot\id\ot\id+\ot\id\ot L\ot\id+\id\ot\id\ot L
$$
and so on. Namely, the matrix $\Delta^k(L)$ satisfies the CH
identity whose roots are $k_1\,\lambda_1+k_2\,\lambda_2$ with
$k_1+k_2=k+1$ and the multiplicity of each root is
$C^{k_1}_{k+1}$. To avoid this multiplicity it suffices to
consider {\it the symmetric component}  (denoted as $L_{(k)}$) of
the matrix $\Delta^{k+1}(L)$ (see below). Finally, the matrix
$L_{(k)}$ has $k+1$ pairwise distinct eigenvalues and its
characteristic polynomial equals to that of $\Delta^k(L)$. Then
by the same method as above we can associate to this matrix $k+1$
idempotents and corresponding  projective modules. Thus, we have
realized the line bundles $E^{k_1,k_2}$ under the guise of
projective modules.

Now, let us pass to the NC variety in question. With matrices
$L_1$ and $L_2$ the commutation relation (\ref{gln}) takes the
form \be L_1\cdot L_2-L_2\cdot L_1=\hbar(L_1P-PL_1) \label{com}
\ee where $P$ is the usual flip. So, we cannot apply the
commutative binomial formula for calculating the powers of
$\Delta^k(L)$. This prevents us from  calculating the CH
identities for the matrices $\Delta^k(L)$ with the above method.

Instead, we will calculate the CH identities directly for {\it
the symmetric components} of these matrices,  also denoted
$L_{(k)},\,k=2,3,...$ and defined as
$$
L_{(k)} = k\,S^{(k)}\,L_1\,S^{(k)},
$$
where $ S^{(k)}$ is the Young symmetrizer in $V^{\ot k}$. Taking
in consideration that the element $S^{(k)}\,L_1$ is already
symmetrized w.r.t. all factors apart from the first one we can
represent the matrix $L_{(k)}$ as
$$
L_{(k)}=S^{(k)}
\,L_1\,(\id+P^{12}+P^{12}P^{23}+...+P^{12}P^{23}...P^{k-1\,k})
$$
where $P^{i\,i+1}$ is the operator transposing the $i$-th and
$(i+1)$-th factors in the tensor product of spaces.

\begin{remark}
We can treat the matrix $L_{(k)}$ as an operator acting from
$V^{\ot k}$ to $V^{\ot k}\ot\Ah$ assuming it to be trivial on all
components except for $V_{(k)}\subset V^{\ot k}$.
\end{remark}

In case $k=2$ we have the following proposition.
\begin{proposition} The CH identity for the matrix $\ll$ restricted to the
symmetric component $V_{(2)}$ of the space $\vv$ is \be
\ll^3-4\hbar\ll^2+4(\alpha+\hbar^2)\ll-8\hbar\alpha\,\id=0.
\label{CHder} \ee
\end{proposition}
\smallskip

{\bf Proof.}\ \ By definition the matrix $\ll$ has the form
$$
\ll = \frac{1}{2}(\id + P^{12})L_1(\id + P^{12}).
$$
Taking in consideration that $L_2=P^{12}\, L_1\,P^{12}$ we
rewrite (\ref{com}) as follows \be L_1P^{12}L_1P^{12} -
P^{12}L_1P^{12}L_1 = \hbar(L_1P^{12} - P^{12}L_1). \label{star}
\ee Now we need some powers of the matrix $\ll$. We find
($L\equiv L_{1}$, $P\equiv P^{12}$):
\begin{eqnarray*}
&&\displaystyle \ll^2 = \frac{1}{2}(\id+P)L(\id+P)L(\id+P)\\
&&\displaystyle \ll^3 =
\frac{1}{2}(\id+P)L(\id+P)L(\id+P)L(\id+P).
\end{eqnarray*}

Then taking into account (\ref{numer}) and (\ref{star}) we have
the following chain of identical transformations for $\ll^3$:
\begin{eqnarray*}
2\ll^3\rlap{$=(\id+P)L\underline{(\id+P)L(\id+P)L}(\id+P) =
(\id+P)L\Bigl[\underline{PLPL}+PL^2+LPL+L^2\Bigr](\id+P)=$}&&\\
&&\hspace*{-18pt} 2 (\id+P)L\Bigl[LPL + \hbar PL - \alpha \id
\Bigr](\id+P) =
4\hbar(\id+P)LPL(\id+P) - 4\alpha (\id+P)L(\id+P)=\\
&&\hspace*{-18pt} 4\hbar (\id+P)L(\id+P)L(\id+P) - 4\hbar
(\id+P)(\hbar L - \alpha\id)
(\id+P) - 4\alpha (\id+P)L(\id+P) = \\
&& \hspace*{-18pt}8\hbar \ll^2 - 8(\hbar^2 +\alpha)\ll +
8\hbar\alpha (\id+P).
\end{eqnarray*}
Canceling the factor $2$ we come to the result \be
\ll^3-4\hbar\ll^2+4(\alpha+\hbar^2)\ll-4\hbar\alpha(\id+P)=0.
\label{prom-res} \ee To complete the proof it remains to note,
that after restriction to the symmetric component $V_{(2)}$ of
the space $\vv$ the last term in (\ref{prom-res}) turns into
$8\hbar\alpha\,\id$ and we come to
(\ref{CHder}).\hfill\rule{6.5pt}{6.5pt}
\medskip

Now, let us  exhibit the matrix $\ll$ in a basis form. In the base
$$
v_{20}=v_1^{\ot 2},\quad v_{11}=v_1\ot v_2+v_2\ot v_1,\quad
v_{02}=v_2^{\ot 2}
$$
of the space $V_{(2)}$ this matrix has the following form \be
\ll=\left(\matrix{2a&2b&0\cr c&a+d&b\cr 0&2c&2d}\right)=
\left(\matrix{2a&2b&0\cr c&0&b\cr 0&2c&-2a}\right)\label{sec} \ee
(the latter equality holds in virtue of the condition $a+d=0$).

In the sequel we prefer to deal with another basis in the space
$V_{(2)}$. Namely, by putting
$$
(v_{20},\,v_{11},\,v_{02})=(u_{20},\,u_{11},\,u_{02})\bf{P}
$$
with transition matrix
$$
\bf{P}=\left(\matrix{0&1&0\cr 1/2&0&-1/2\cr-i/2&0&-i/2}\right)
$$
we transform $\ll$ into the form \be
\oll={\bf{P}}{\ll}{\bf{P}^{-1}} ={2}\left(\matrix{0&-z&y\cr
z&0&-x\cr -y&x& 0}\right)\label{lov}. \ee

The matrix $\oll$ is expressed through the generators $(x,y,z)$
and it is better adapted to the compact form of the NC varieties
in question. However, it satisfies the same NC version of CH
identity (\ref{numer}).

\begin{remark} By straightforward checking we can see that
the roots of (\ref{CHder}) are \be \lambda_{20}=\h -
\sqrt{\h^2-4\alpha}=2\lambda_1\,\,\lambda_{11}=2\h=2(\lambda_1+\lambda_2)
,\,\,{\rm and} \,\,\lambda_{02}=\h
+\sqrt{\h^2-4\alpha}=2\lambda_2.\label{three} \ee These
quantities are "eigenvalues" of the matrix (\ref{sec}) which is
the restriction of the matrix (\ref{first}) to the symmetric part
of the space $\vv$. A similar restriction of the matrix
(\ref{first}) to the skewsymmetric part of the space $\vv$ gives
rise to the operator
$$
v_1\ot v_2-v_2\ot v_1\to (v_1\ot v_2-v_2\ot v_1)\ot (a+d)=0.
$$
Thus, the matrix (\ref{first})  on the whole space $V^{\ot 2}$
has 4 distinct eigenvalues (those (\ref{three}) and 0) and in
contrast with the commutative case its minimal polynomial cannot
be of the third degree.
\end{remark}

As for higher extensions $L_{(k)},\,k>2$ of the matrix $L$
 the following holds.

\begin{proposition}\label{x}
For any integer $k>2$ there exists a  polynomial
$$P_k(x)=\lambda^k+\sum_{i=1}^{k}a_{k-i}\,\lambda^{k-i}$$
with numerical coefficients such that $P_k(L_{(k)})=0$. Moreover,
its roots are
$$\lambda_{k_1k_2}=k_1\lambda_1+k_1\,k_2\,(\lambda_1+\lambda_2)
+k_2\lambda_2,\,\,\, k_1+k_2=k.$$
\end{proposition}

This formula can be deduced from \cite{Ro}. We will present its
q-analogue in \cite{GLS}.

Similar to the basic case the roots of the polynomial $P_k$ will
be called eigenvalues of the corresponding matrix $L_{(k)}$.

Now, let us assume the eigenvalues $\lambda_{20}, \lambda_{11},
\lambda_{02}$ to be also pairwise distinct. Then, by using the
same method as above we can introduce the idempotent
$$
e_{20}=(\lambda_{11}\id-
L_{(2)})(\lambda_{02}\id-L_{(2)})/(\lambda_{11}-
\lambda_{20})(\lambda_{02}-\lambda_{20})
$$
and similarly $e_{11}$ and $e_{02}$ corresponding to the
eigenvalues $\lambda_{11}$, and $\lambda_{02}$ respectively. The
related q.l.b. (projective $\Ah$-modules) will be denoted
$E^{20}(\h)$, $E^{11}(\h)$, and $E^{02}(\h)$ respectively.

Assuming the eigenvalues of the polynomials $P_k, \,k>2$ (see
proposition \ref{x}) to be also pairwise distinct we can
associate to the matrix $L_{(k)}$  $k+1$ idempotents
$$
e_{k_1k_2},\,\,k_1+k_2=k,\,\,k>0
$$
and the corresponding q.l.b. $E^{k_1,k_2}(\h)$. For $k=0$ we set
$e_{00}=1$. The corresponding q.l.b. is $E^{0,0}(\h)$.

\begin{remark}
Let us note that if we do not fix any value of $\triangle$ we can
treat the elements $\tr\,e_{ij}$ as those of
$M_2(U(sl(2)_{\h}))\ot R$ where $R$ is the field of fractions of
the algebraic closure $\overline{Z[U(sl(2)_{\h})]}$.
\end{remark}

\begin{remark}
Note that there exists a natural generalization of the above
constructions giving rise to some "braided varieties" and
corresponding "line bundles" as follows. Let  $R$ be a  Hecke
symmetry of rank 2 (cf. \cite{G}).
 Then the matrix $L$ satisfying the
RE with such $R$ obeys an equation analogous to (\ref{CH}) but
with appropriate trace and determinant (cf. \cite{GPS}).
Introducing the quotient algebra $\Ah$ in a similar way we treat
it as a braided analogue of a NC sphere. Then, by defining the
extensions $L_{(k)}$ as above (but with modified meaning of the
symmetric powers of the space $V$) we can define a family of
q.l.b. over such a "NC braided variety" as above. This
construction will be presented in details in \cite{GLS}.
\end{remark}

Now, we pass to computing the quantities $\tr \,e_{k_1k_2}$.

\begin{proposition} The following relation holds
$$
\tr \,e_{k_1k_2}=1+{{(k_1-k_2)\h}\over{\sqrt{\h^2-4\alpha}}}.
$$
\end{proposition}
A proof of this formula will be given in \cite{GLS} in a more
general context including its q-analogue.

\section{Isomorphic modules and multiplicative structure}
\setcounter{equation}{0}

First of all we discuss the problem of isomorphism between the
projective modules introduced above (namely, q.l.b.). There
exists a number of definitions of isomorphic modules over ${\Bbb
C}^*$-algebras (cf. \cite{W}). However, for the algebras in
question we use the following definition motivated by \cite{R}.

\begin{definition} We say that two projective modules
$M_1\subset\A^{\oplus m}$ and $M_2\subset\A^{\oplus n}$ over an
algebra $\A$ corresponding to the idempotents $e_1$ and $e_2$
respectively are isomorphic iff there exist two matrices $A\in
M_{m,n}(\A)$ and $B\in M_{n,m}(\A)$ such that
$$
AB=e_1,\quad BA=e_2,\quad A=e_1 A=Ae_2,\quad B=e_2 B=B e_1.
$$
\end{definition}

\begin{proposition}
The q.l.b. $E^{1,1}(\h)$ is isomorphic to that $E^{0,0}(\h)$.
\end{proposition}
{\bf Proof}\ \ It is not difficult to see that
$$e_{11}=(L_{(2)}^2-2\h L_{(2)}+4\alpha\id)/(4\alpha).$$

Then,  by straightforward calculations we check that for the
idempotent $e_{11}$ the following relation holds
$$
(4\alpha)^{-1}\left(\left(\matrix{4a^2+2bc&4ab&2b^2\cr
2ca&2cb+2bc&-2ba\cr
2c^2&-4ac&2cb+4a^2}\right)-2\h\left(\matrix{2a&2b&0\cr c&0&b\cr
0&2c&-2a }\right)+ 4\alpha\left(\matrix{1&0&0\cr 0&1&0 \cr
0&0&1}\right)\right)=
$$
$$
(4\alpha)^{-1} \left(\matrix{-2b\cr 2a\cr
2c}\right)\left(\matrix{c&-2a&-b}\right).
$$

Passing to the matrix $\oll$ we get
$$
e_{11}=(4\alpha)^{-1}(\oll^2-2\h\oll+4\alpha\id)=\alpha^{-1}
\left(\matrix{x\cr y\cr z}\right)\left(\matrix{x&y&z}\right)
$$
(So, the idempotent $e_{11}$ defines the following operator in
${\A}^{\oplus 3}$:
$$\left(\matrix{f_1\cr f_2\cr f_3}\right)\mapsto
(\alpha)^{-1} \left(\matrix{x\cr y\cr z}\right)(x\cdot f_1+y\cdot
f_2 +z\cdot f_3).)$$

It remains to say that if we put
$A=(\alpha)^{-1}\left(\matrix{x&y&z}\right)$ and
$B=\left(\matrix{x\cr y\cr z}\right)$ we satisfies the definition
above with $e_1=e_{00}$ and $e_2=e_{11}$ \hfill\rule{6.5pt}{6.5pt}

In general, the problem of isomorphism between modules
(\ref{isom}) is open. We can only conjecture that the q.l.b.
(\ref{isom}) are isomorphic to each other.

\begin{remark}
If the algebra $\A$ is not commutative the first two relations of
the definition do not yield the equality $\tr\, e_1=\tr\, e_2$.
However, proposition 8 implies that
$$
\tr \,e_{k_1,k_2}=\tr\,e_{k_1+l,k_2+l}.
$$
\end{remark}

Now, we can introduce an associative product on the set of q.l.b.
in a natural way by setting
$$
E^{k_1,k_2}(\h)\cdot E^{l_1,l_2}(\h)=E^{k_1+l_1,k_2+l_2}(\h).
$$
This product is evidently associative and commutative. In
particular, we have
$$E_1(\h)\cdot E_2(\h)=E^{1,1}(\h),\,\,E_1(\h)\cdot E_1(\h)=E^{2,0}(\h),
\,\, E_2(\h)\cdot E_2(\h)=E^{0,2}(\h).$$

The family of the modules $E^{k_1,k_2}(\h)$ equipped with this
product is denoted $\prePic(\Ah)$ and called {\em prePicard
semigroup} of the NC sphere.

Assuming that the  q.l.b. (\ref{isom}) are indeed isomorphic to
each other we can naturally define the Picard group $\Pic(\Ah)$
of the NC sphere as the classes of isomorphic modules
$E^{k_1,k_2}(\h)$ equipped with the above product.

So, under this assumption, the Picard group $\Pic(\Ah)$ of the NC
sphere is at most $\Z$ (recall that $\Pic(S^2)=\Z$).

Now, consider the problem of computing the pairing \be
<\,\,,\,\,>\, :\, \prePic(\Ah)\ot K^0\to\Bbb K. \label{coupl} \ee
Such a pairing plays the key role in the Connes version of the
index formula. (Usually, one considers $K_0$ instead of $\prePic$
but we restrict ourselves to "quantum line bundles". Moreover,
assuming the conjecture formulated before remark 12  to be true
we can replace $\prePic$ in formula (\ref{coupl})  by $\Pic$.) Let
us remind that $K^0$ stands for the Grothendieck ring of the
category of irreducible modules of the algebra in
question\footnote{By this we mean the Grothendieck ring of the
algebra $U(su(2)_{\h})$ or what is the same $U(su(2))$ since $\h$
does not matter here. However, any irreducible
$U(su(2)_{\h})$-module defines a relation between $\h$ and
$\alpha$ (see (\ref{alp})) and therefore the factors in the
formula (\ref{index}) are not independent.}.  In the spirit of
the NC index (cf. \cite{L}) the paring (\ref{coupl}) can be
defined as \be <E^{k_1,k_2}(\h), U>=\tr
\,\pi_U(\tr(e_{k_1,k_2}))\label{index} \ee where
$\tr(e_{k_1,k_2})\in\Ah$ and $\pi_U:\Ah\to\End(U)$ is the
representation corresponding to the irreducible $U$.

It is not difficult to see that the result of the  pairing of the
module $E^{k_1,k_2}(\h)$ with the irreducible  $U_j$ of the spin
$j$  is equal to the quantity $n\,\tr\,e_{k_1,k_2}$ evaluated at
the point \be \alpha=-\h^2(n^2-1)/4\quad {\rm where}\quad
n=\dim\,U_j=2j+1.\label{alp} \ee

In particular, we have
$$
<E^{0,0}(\h), U_j>=n,\,\,<E^{1,0}(\h), U_j>=n+1,\,\,<E^{0,1}(\h),
U_j>=n-1,
$$
$$
<E^{2,0}(\h), U_j>=n+2,\,\, <E^{0,2}(\h), U_j>=n-2
$$
(here we assume that $\sqrt{n^2\h^2}=n\h$). More generally, if
$n>k_1+k_2$ we have
$$<E^{k_1,k_2}(\h), U_j>=k_1-k_2+n.$$
This formula follows immediately from proposition 9. There exists
a q-analogue of this formula which will be considered in
\cite{GLS}.

\section{Differential calculus via projective modules}
\setcounter{equation}{0}

The aim of this section is to develop some elements of
differential calculus on the NC sphere in terms of projective
modules (namely, the q.l.b. above and their direct sums).
Usually, a differential algebra associated to a NC algebra is
much bigger than the classical one even if such an algebra is a
deformation of the commutative coordinate ring of a given variety
(cf. for example \cite{GVF}). Moreover, the components of these
differential algebras are not finitely generated modules. This
leads to a NC version of de Rham complex which is drastically
different from the classical one.

We suggest another  NC version of de Rham complex looking like
its classical counterpart.  The main idea of our approach
consists in the following. Instead of using the Leibniz rule and
introducing a way of transposing "functions" and "differential"
we treat any term of the classical de Rham complex as a
projective module. Let us consider projective $\Ah$-modules which
are analogues of the former ones.
 By decomposing them into direct sums of irreducible $\Ah$-modules
we define a differential on each irreducible component following
the classical pattern (a detailed description of this
decomposition in the classical case is given in \cite{AG}). Then
the property $d^2=0$
 for the deformed differential is preserved automatically
and the cohomology of the final complex are just the same as in
the classical case. Let us describe this approach in more details.

However,  we want to begin with a NC analogue of the tangent
bundle $T(S^2)$ on the sphere. Since this bundle is complementary
to the normal one and since the latter bundle (treated as a
module) is nothing but $E^{1,1}$, it is natural to define NC
analogue of $T(S^2)$ as
$$
T(\Ah)=E^{2,0}(\h)+E^{0,2}(\h).
$$
We call it {\em a tangent module} on the NC sphere.

This module can be represented by the equation $\Im\, e_{11}=0$.
It is equivalent to the relation \be
u_{20}\,x+u_{11}\,y+u_{02}\,z=0. \label{rel} \ee This means that
the module $T(\Ah)$ is realized as the quotient of the free
module $\Ah^{\oplus 3}$ generated by the elements \be
u_{20},\,u_{11},\,u_{02} \label{gener} \ee over the submodule
$\{Cf,\, f\in\Ah\}$ where $C$ is the l.h.s. of (\ref{rel}).

In the classical case relation (\ref{rel}) is motivated by an
operator meaning of the tangent space. Namely, if generators
(\ref{gener}) are treated as infinitesimal rotations of the
sphere\footnote{By means of the Kirillov bracket we can represent
these rotations as $u_{20}=\{x,\,\cdot\,\}$ and so on.} and the
symbols $x,\, y,\, z$  in  (\ref{rel}) are considered as
operators of multiplication on the corresponding functions, then
the element $C$ treated as an operator is trivial. This allows us
to equip the tangent module $T(S^2)$ with an operator meaning by
converting any element of this module into a vector field. Thus,
we get the action
$$
\A\ot T(S^2)\to \A,\,\, \A={\Bbb K}(S^2)
$$
which consists in applying a vector field to a function.

However, if we assign the same meaning to the generators
(\ref{gener}) and to those $x,\, y,\,z$  on the NC sphere (by
setting $u_{20}=[x,\,\cdot\,]$ and so on and considering $x,\,
y,\,z$ as operators of multiplication on the corresponding
generator) then the element $C$ treated as an operator  is no
longer trivial. This is the reason why we are not able to provide
the tangent module $T(\Ah)$ with a similar action on the algebra
$\Ah$.

\begin{remark}
We want to stress that the space of derivations of the algebra
$\Ah$ often considered as a proper NC counterpart of the usual
tangent space does not have any $\Ah$-module structure. So, for a
NC variety which is a deformation of a classical one we have a
choice: which properties of the classical object we want to
preserve. In our approach we prefer to keep the property of the
tangent space to be a projective module. In the same manner we
will treat the terms of the deformed de Rham complex (see below).
\end{remark}

Passing to the cotangent bundle $T^*(S^2)$ or otherwise stated to
the space of the first order differentials $\Om^1(S^2)$ we see
that it is isomorphic to $T(S^2)$. Therefore, it is defined by
the same relation (\ref{rel}) but with another meaning of
generators (\ref{gener}): now we treat them as the differentials
of the functions $x,\, y$ and $z$ respectively:
$$
u_{20}=d\, x,\,u_{11}=d\, y,\,u_{02}=d\, z.\label{mean}
$$

This gives us a motivation to introduce {\em cotangent module}
$T^*(\Ah)$ on the NC sphere similarly to the tangent one but with
a new meaning of generators (\ref{gener}). We will also use the
notation $\Om^1(\Ah)$ for the module $T^*(\Ah)$  and call the
elements of this space the first order differential forms on the
NC sphere.

Up to now we have no modifications in the relations defining the
tangent and cotangent objects. Nevertheless, it is no longer so
for the NC counterpart of the second order differentials space
$\Om^2(S^2)$. In the classical case this space is defined by \be
u_{11}\,z-u_{02}\, y=0,\quad -u_{20}\,z+u_{02}\, x=0,\quad
u_{20}\,y-u_{11}\, x=0 \label{syst} \ee with the following
meaning of generators (\ref{gener}) \be u_{20}=d\, y\wedge d\, z,
\,u_{11}=d\, z\wedge d\, x,\,u_{02}=d\, x \wedge d\, y.
\label{diff} \ee In a concise  form we can rewrite (\ref{syst})
as \be (u_{20},\,u_{11},\,u_{02})\cdot \oll=0 \label{newrel} \ee
with $\oll$ given by (\ref{lov}).

However, in the NC case we should replace (\ref{newrel}) by
$$
(u_{20},\,u_{11}\,u_{02})\cdot \oll=2\h (u_{20},\,u_{11}\,u_{02})
$$
or in more detailed form
$$
u_{11}\,z-u_{02}\, y=\h\, u_{20},\quad -u_{20}\,z+u_{02}\, x=\h\,
u_{11},\quad u_{20}\,y-u_{11}\, x=\h \,u_{02}.
$$
This is motivated by the fact that $2\h$ becomes  an eigenvalue
of the matrix $\oll$. Let us denote the corresponding quotient
module $\Om^2(\Ah)$. Of course, we can represent generators
(\ref{gener}) in form (\ref{diff}) but now it does not have any
sense. Nevertheless, we will call the elements of this module
{\it the second order differential forms}.

Thus, we have defined all terms of the following {\em de Rham
complex} on the NC sphere \be 0\to
\Om^0(\Ah)=\Ah\to\Om^1(\Ah)\to\Om^2(\Ah)\to 0. \label{Rham} \ee
These terms are projective $\Ah$-modules which are deformations
of the corresponding $\A$-modules. Now we would like to define a
differential $d$ in this complex in such a way that this complex
would have just the same cohomology as in the classical case.

Let us decompose each term of the classical de Rham complex into
a direct sum of  irreducible  $SU(2)$-modules (cf. \cite{AG}).
Since in the classical case the differential $d$ is
$SU(2)$-covariant then any such an irreducible module is mapped
by the differential either to 0 or to an isomorphic module.

Remark that  the terms of (\ref{Rham}) being a deformation of
their classical counterparts  consist of just the same irreducible
$SU(2)$-modules as in the classical case. This property
 follows from the fact that the
modules in question are projective and the corresponding
idempotents smoothly depend on $\h$.

Now, we define the differential in (\ref{Rham}) as a mapping
which takes any irreducible $SU(2)$-module to 0 or to the
isomorphic module following the classical pattern.

As an example we  describe how the differential $d$ acts on the
algebra $\Ah$ itself. It is easier to use the non-compact form of
the algebra. The algebra $\Ah$ is a multiplicity free direct sum
of irreducible $SL(2)$-modules. Their highest (or lowest) weight
elements w.r.t. the action of the group $SL(2)$ are $b^k,\,
k=1,2...$. In virtue of $SL(2)$-covariance it suffices to define
the differential on these elements. Similar to the classical
case  we set
$$
d\,b^k=k\,(d\,b)\, b^{k-1}.
$$

Stress that this relation is not obtained as a result of
transposing "functions" and "differentials" mixed up in virtue of
the Leibniz rule but it is imposed by definition.  In our version
of NC de Rham complex
 we do not use either any form of the Leibniz rule or any
transposing "functions" and "differentials".

However, by construction we have $d^2=0$ and just the same
cohomology as in the classical case. Similarly to the classical
case this cohomology is generated by 1 in the  term $\Om^0(\Ah)$
and by the element $u_{20}\, x+u_{11}\, y+u_{20}\, z$ in the term
$\Omega^2(\Ah)$.


\begin{thebibliography}{DGMMM}

\bibitem[AG]{AG} P.Akueson, D.Gurevich {\em Some aspects of braided
geometry: differential calculus, tangent space, gauge theory}, J
Phys. A: Math. Gen 32 (1999), pp 4183--4197.

\bibitem[C]{C} A.Connes  {\em Noncommutative Geometry}, Academic Press,
1994

\bibitem[DL]{DL} L.Dabrowski, G.Landi {\em Instanton algebras and
quantum 4 spheres}, q-alg/0101177.

\bibitem[G-T]{G-T} I.Gelfand, D.Krob, A.Lascoux, B.Leclerc,
V.Retakh, J-Y.Thibon {\em Noncommutative symmetric functions},
Adv. Math. 112 (1995), pp.218--348.

\bibitem[Go]{Go} M.D.Gould {\em On the matrix elements of the
$U(n)$ generators}, J.Math.Phys. 22 (1981), pp.15--22.

\bibitem[GVF]{GVF} J. Gracia-Bondia, J.Varilly, H.Figueroa
{\em Elements of Noncommutative Geometry}, Birkh\"auser, 2000.

\bibitem[G]{G} D.Gurevich {\em Algebraic aspects of the quantum
Yang-Baxter equation}, Leningrad Math.J. 2 (1991), pp. 801--828.

\bibitem[GLS]{GLS} D.Gurevich, R.Leclercq, P.Saponov {\em
Equivariant noncommutative index on braided sphere}, in
preparation.

\bibitem[GS]{GS} D.Gurevich, P.Saponov {\em Quantum line bundles via
Cayley-Hamilton identity}, J. Phys. A : Math. Gen.  34 (2001),
pp. 4553 -- 4569.

\bibitem[GPS]{GPS} D.Gurevich, P.Pyatov, P.Saponov {\em Hecke
Symmetries and Characteristic Relations on Reflection Equation
Algebra}, Lett. Mat. Phys. 41 (1997) pp. 255--264.

\bibitem[H]{H} T.Hodes {\em Morita isomorphism of Primitive factors of
$U(sl(2))$}, Contemporary Mathematics 139 (1992), pp.175--179.

\bibitem[Ha]{Ha} P.Hajac {\em Bundles over the quantum sphere
and noncommutative index theorem}, K-theory 21 (2001),
pp.141--150.

\bibitem[KT]{KT} I.Kantor, I.Trishin {\em On concept of Determinant
in the Supercase}, Comm. in Algebra 22 (1994), pp. 3679--3739.

\bibitem[K]{K} A.Kirillov {\em Introduction to family algebras},
Moscow Math.J. 1 (2001) pp.49--64.

\bibitem[KS]{KS} A.Konechny, A.Schwarz {\em Introduction to M(atrix)
theory and noncommutative geometry}, hep-th/0012145.

\bibitem[LM]{LM} G.Landi, J.Madore {\em Twisted configurations
over Quantum Euclidean Spheres}, preprint q-alg/0102195.

\bibitem[L]{L} J-L.Loday {\em Cyclic homology}, Springer, 1998.

\bibitem[NT]{NT} M.Nazarov, V.Tarasov {\em Yangians and Gelfand-Zetlin
bases}, Publ. RIMS 30 (1994) pp.459--478.

\bibitem[Ri]{Ri} M.Rieffel {\em Projective modules over higher-dimensional
non-commutative tori}, Can.J. Math. 40 (1988) pp.257--338.

\bibitem[Ro]{Ro} N.Rozhkovskaya {\em Family algebras of representations
with simple spectrum}, preprint ESI-Vienna 1045 (2001).

\bibitem[R]{R} J.Rosenberg {\em Algebraic K-tehory and its applications},
Springer, 1994.

\bibitem[Se]{Se} J.-P.Serre {\em Modules projectifs et espaces fibr\'es a
fibre vectorielle}, Seminaire Dubreil-Pisot, Fasc.2, Expos\'e 23
(1957/1958).

\bibitem[Sw]{Sw} R.Swan {\em Vector bundles and projective modules}, Trans.
AMS 105 (1962), pp. 264--277.


\bibitem[W]{W} N.Wegge-Olsen {\em K-theory and ${\Bbb C}^*$-Algebras},
Oxford University Press, 1993.


\end{thebibliography}
\end{document}